\documentclass[12pt,a4paper]{article}
\usepackage{euler,natbib,amsmath,defns}
\allowdisplaybreaks
\IfFileExists{ajr.sty}{\usepackage{ajr}}{}

\newcommand{\E}[2]{$#1\cdot10^{#2}$}
\newcommand{\pbc}{\textsc{tbc}}

\begin{document}
    
\title{General tooth boundary conditions for equation free modelling}

\author{A.~J. Roberts\thanks{Computational Engineering and Sciences
Research Centre, Department of Maths \& Computing,
University of Southern Queensland, Toowoomba, Queensland~4352,
\textsc{Australia}.  \protect\url{mailto:aroberts@usq.edu.au}} \and
I.~G. Kevrekidis\thanks{Program in Applied and Computational
Mathematics, Princeton University, Princeton, NJ~08544, USA.
\protect\url{mailto:yannis@Princeton.edu}}}

\maketitle

\begin{abstract}
We are developing a framework for multiscale computation which enables
models at a ``microscopic'' level of description, for example Lattice
Boltzmann, Monte Carlo or Molecular Dynamics simulators, to perform
modelling tasks at ``macroscopic'' length scales of interest.
The plan is to use the microscopic rules restricted to small ``patches" of
the domain, the ``teeth'', using interpolation to bridge the ``gaps".
Here we explore general boundary conditions coupling the 
widely separated ``teeth'' of the microscopic simulation that
achieve high order accuracy over the macroscale.
We present the simplest case when the microscopic simulator is 
the quintessential example of a partial differential equation.
We argue that classic high-order interpolation of the macroscopic field
provides the correct forcing in whatever boundary condition 
is required by the microsimulator.
Such interpolation leads to Tooth Boundary Conditions which achieve
arbitrarily high-order consistency.
The high-order consistency is demonstrated on a class of linear partial
differential equations in two ways: firstly through the eigenvalues of
the scheme for selected numerical problems; and secondly using the dynamical systems
approach of holistic discretisation on a general class of linear
\textsc{pde}s.
Analytic modelling shows that, for a wide class
of microscopic systems, the subgrid fields and the effective
macroscopic model are largely independent of the tooth size and the
particular tooth boundary conditions.
When applied to patches of microscopic simulations these tooth boundary
conditions promise efficient macroscale simulation.
We expect the same approach will also accurately couple patch
simulations in higher spatial dimensions.
\end{abstract}

\paragraph{Keywords:} multiscale computation, gap tooth scheme,
coupling boundary conditions, high order consistency

\tableofcontents

\section{Introduction}

The components of physical systems often operate on vastly different
space and time scales~\cite*[]{Dolbow04}.
We must somehow simulate such systems on the scale of interest and
operation.
But systems that depend on physical processes at multiple scales pose
notorious difficulties.
These multiscale difficulties are major obstacles to
progress in fields as diverse as environmental and geosciences, climate,
materials, combustion, high energy density physics, fusion, bioscience,
chemistry, power grids and information networks~\cite[]{Dolbow04}.

Here we further develop the equation free approach to
multiscale modelling~\cite[]{Kevrekidis03b}.
Given a numerical simulator for physical components at much smaller
scales than the scale of primary interest, the aim of the methodology
is to bridge the space and time scales to simulations resolving the
macroscale of interest.
Here we focus on bridging \emph{space scales} by improving the accuracy of
the gap-tooth methodology for microsimulators~\cite[]{Gear03,
Samaey03a, Samaey03b}.
Crucially, our gap-tooth methods must adapt to whatever microsimulator
code is provided; one key application of this work is to
microsimulators that are tried and tested legacy codes that we do not want 
to modify.

The equation-free approach provides \emph{on the fly} 
closure methods which constitute critical components of, for
example, mathematical homogenization~\cite[e.g.]{Samaey03b,
Gustafsson03, Balakotaiah03}, renormalization group
techniques~\cite[e.g.]{Ei99, Mudavanhu03, Chorin05}, and multiscale
finite elements~\cite[e.g.]{Hou97, Chen02}.
These closure methods not only need to be computationally efficient but
also need to be capable of reproducing the physical dynamics with high
fidelity.
That is, we seek a methodology that can be systematically refined.

\begin{figure}
    \centering
    \includegraphics{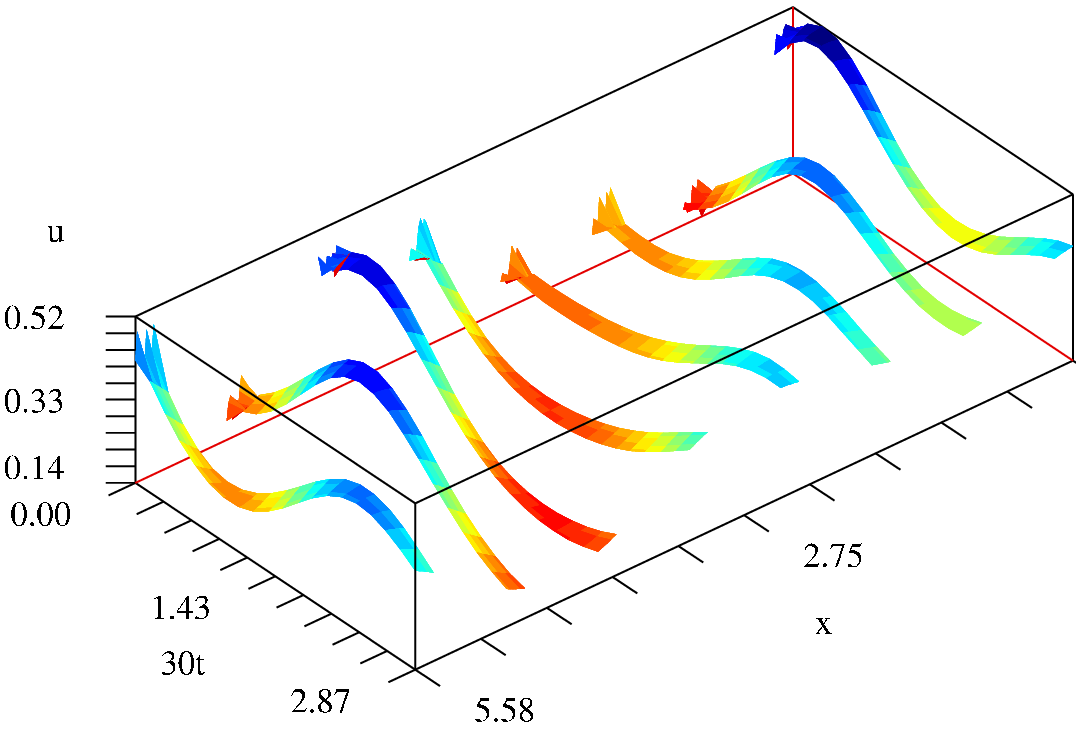}
    \caption{simulation of Burgers' equation using Dirichlet boundary
    conditions on the teeth (specified~$u$ on the edges).}
    \label{fig:burg3}
\end{figure}

Using microscopic simulators of the one dimensional Burgers' equation,
\cite{Roberts04d} demonstrated the possibility of achieving high order
accuracy in the gap-tooth scheme for macroscale dynamics.
The particular microsimulator we use is a fine scale discretization of the
\pde\ which we execute only in the interior of the teeth 
(see Figure~\ref{fig:burg3}).
At each time step during execution, the microsimulator within each tooth requires
boundary values which must be continuously updated.
If the microsimulator
was to be executed over the entire macrodomain, these boundary values would
naturally come from the immediately neighboring fine grid; this grid
is missing in gap-tooth simulation.
That pilot study only considered microsimulators which had boundary
conditions of specified \emph{flux} at the edges of their simulation teeth.
Here we generalise the analysis to consider microsimulators with either
\begin{itemize}
	\item Dirichlet boundary conditions of specified field~$u$ at the
	tooth edges, Section~\ref{sec:diri},

	\item mixed boundary conditions of specified $av_j\pm
	b\partial_xv_j$ at the tooth edges, Section~\ref{sec:mixed}, or

	\item nonlocal two-point boundary conditions such as those arising in a
	microscale discretisation of a \pde, Section~\ref{sec:twopt}.
\end{itemize}

Consider the gap-tooth scheme \cite[e.g.]{Gear03, Samaey03a}
illustrated in Figure~\ref{fig:burg3}.
Let $v_j(x,t)$~be the fine scale, microscopic field in the $j$th~tooth, 
and $U_j$~the $j$th~coarse grid value; that is, the value at the center
of each tooth.
Let the tooth width be~$h$.
Then the edge of a tooth lies at a distance~$h/2$ from its coarse grid point, a
fraction $r=h/(2H)$ to the next coarse grid point.
The amount of computation performed by microsimulators is proportional
to the width of the (microscale) teeth.
Hence we aim for the fraction~$r$ to be as small as possible,
so that the teeth are a
relatively small part of the physical domain and the computational cost
minimised.
The coupling rule developed in Sections~\ref{sec:diri},
\ref{sec:mixed}~and~\ref{sec:twopt} is that you \emph{obtain whatever
values are necessary for the boundaries of the microscopic simulators by classic
interpolation of the \emph{macroscopic} grid values from neighbouring teeth.}
As a nonlinear example of the coupling we develop, Figure~\ref{fig:burg3}
shows a gap-tooth simulation of the nonlinear dynamics of Burgers' equation in one
spatial dimension.

This coupling rule promotes a strong connection between classic finite
difference discretisations of \pde{}s, classic finite elements, and the
methodology of the gap-tooth scheme.
First, the \pde\ acts as the quintessential example of a microsimulator
in that it informs us of the dynamics in an `infinitesimal patch'.
The only difference between the gap-tooth scheme and the spatial
discretisation of \pde{}s is that the microsimulators in the gap-tooth
scheme encode the dynamics on \emph{small finite patches}, whereas the \pde\
encodes the dynamics on \emph{infinitesimal} patches.
Consequently, classic interpolation serves the same role in both:
namely, the interpolation appropriately transfers information from the
macroscale of interest to the microscale simulators.
Second, the theoretical support for the gap-tooth scheme is based upon
a subgrid scale structure, as is the classic finite element method.
Also, the solvability condition in the construction of the theoretical
gap-tooth model is similar to the Galerkin projection of finite elements.
But whereas finite elements \emph{impose} a class of subgrid fields, both the
theoretical approach here and the gap-tooth scheme use actual subgrid
scale dynamics, obtained from the microsimulator or the \pde, to \emph{obtain}
appropriate subgrid scale fields.
Thus this approach systematically implements model closures for macroscale
discretisations.

In Section~\ref{sec:ind} we prove that classic interpolation connects accurately
the teeth across the gaps for the general \emph{linear} fourth order \pde.
The technique of holistic discretisation \cite[e.g.]{Roberts01g}
resolves subgrid scale structures to reproduce with high fidelity the
dynamics of specified \pde{}s \cite[]{Roberts00a}.
The techniques were adapted by \cite{Roberts04d} to the gap-tooth
scheme in the case where the microsimulator requires Neumann boundary
conditions of specified slope\slash flux (see Section~\ref{sec:mixed}).
Using the same techniques, Section~\ref{sec:ind} analyses the use of classic
interpolation of macroscale grid values in the microscale simulators
and shows the following desirable properties:
\begin{itemize}
    \item the approach generates macroscopic discretisations which are
    consistent with the microscopic dynamics to high order in the
    macroscopic tooth separation~$H$;
    
	\item the macroscopic model and the microscopic solution field are
	essentially independent of the size of the teeth, measured by~$r$;
	and

	\item the macroscopic model and the microscopic solution field are
	essentially independent of the details of the tooth boundary
    conditions (\pbc{}s) that couple the teeth together.
\end{itemize}
Thus our proposed rule generates gap-tooth schemes that may be
systematically refined to high order accuracy, and gives rise to macroscale
simulations that are largely independent of irrelevant microscale
parameters.

\section{Dirichlet teeth (specified~$u$)}
\label{sec:diri}

In this section we consider the case of microsimulators that require
at each time step
the field values on the edge of each spatial patch to be specified.
We model this case by \pde{}s with Dirichlet conditions coupling the
dynamics in the teeth.
The values for these Dirichlet conditions are obtained by 
interpolation across the gaps between the teeth using finite
difference operators and the exact relationships between the operators.
When applied to simple diffusion, the resulting scheme has
high order accuracy.

Discrete operators are essential in the analysis.  
Define the shift
operator $Eu(x)=u(x+H)$ and equivalently $EU_j=U_{j+1}$ as appropriate
for steps on the coarse grid size~$H$.  
Then we use the following
identities for discrete operators \cite[p.65, e.g.]{npl61}:
\begin{eqnarray}
    \text{mean}&& \mu = \half(E^{1/2}+E^{-1/2})\,,\\
    \text{difference}&& \delta = E^{1/2}-E^{-1/2}\,,\\
    \text{shift}&& E = 1+\mu\delta+\half\delta^2\,,\\
    \text{derivative}&& H\partial_x = 2\sinh^{-1}\half\delta
    = \delta-\rat16\delta^3+\Ord{\delta^5}\,, \label{eq:derdel}\\
    &&\mu^2 = 1+\rat14\delta^2\,.\label{eq:mudel}
\end{eqnarray}
%
Formulae involving these operators become more accurate as the
differences~$\delta$ become small.  
Such small differences arise either
as the macroscopic grid size $H\to0$ or equivalently as the gradients
of the physical field~$u$ become small.  

For example, \cite{Roberts04d}
showed arbitrary order consistent macroscopic dynamics from a gap-tooth
scheme as the grid size $H\to0$\,.
The key to that analysis is the following transformation of the operator
for evaluating spatial derivatives~$H\partial_x$ at the patch
boundaries~$E^{\pm r}$:
\begin{displaymath}
    E^{\pm r}H\partial_x 
    =(1+\mu\delta+\half\delta^2)^{\pm r}2\sinh^{-1}\half\delta\,.
\end{displaymath}
But this right-hand side, when expanded in a Taylor series in small
differences~$\delta$, is composed of terms which have an odd number of
centred operators $\delta$~and~$\mu$.
Consequently the right-hand side above would require field values
halfway between the grid values.
These are unknown.
Instead, from~\eqref{eq:mudel}, multiply the right-hand side by the
identity~$\mu/\sqrt{1+\delta^2/4}$, and then expand in small
differences~$\delta$:
\begin{eqnarray}
    E^{\pm r}H\partial_x 
    &=&(1+\mu\delta+\half\delta^2)^{\pm r}2\sinh^{-1}\half\delta
    \nonumber\\
    &=&\frac\mu{\sqrt{1+\rat14\delta^2}}
    (1+\mu\delta+\half\delta^2)^{\pm r}2\sinh^{-1}\half\delta
    \nonumber\\
    &=&\mu\delta \pm r\delta^2 
    -(\rat16-\rat12r^2)\mu\delta^3
    \mp r(\rat1{12}-\rat16r^2)\delta^4
    \nonumber\\&&{}
    +(\rat1{30}-\rat18r^2+\rat1{24}r^4)\mu\delta^5
    \pm r(\rat1{90}-\rat1{36}r^2+\rat1{120}r^4)\delta^6
    \nonumber\\&&{}
    -(\rat1{140} -\rat7{240}r^2 +\rat1{72}r^4 -\rat1{720}r^6 )\mu\delta^7
    \nonumber\\&&{}
    \mp r(\rat1{560} -\rat7{1440}r^2 +\rat1{480}r^4 -\rat1{5040}r^6 )\delta^8
    +\Ord{\delta^9}.
    \label{eq:boper}
\end{eqnarray}

For microsimulators with Dirichlet boundary conditions, we adapt the
earlier analysis of \cite{Roberts04d}.  
But instead of determining the slopes at the tooth boundaries as above, 
the following interpolation of the
macroscopic field determines the field values~$u$ on the edges of the
teeth:
\begin{eqnarray}
    E^{\pm r}
    &=&(1+\mu\delta+\half\delta^2)^{\pm r}
    \nonumber\\
    &=&1 \pm r\mu\delta +\rat12r^2\delta^2 
    \pm \rat1{3!}r(r^2-1)\mu\delta^3
    +\rat1{4!} r^2(r^2-1)\delta^4
    \nonumber\\&&{}
    \pm\rat1{5!}r(r^2-1)(r^2-4)\mu\delta^5
    +\rat1{6!}r^2(r^2-1)(r^2-4)\delta^6
    \nonumber\\&&{}
    \pm \rat1{7!}r(r^2-1)(r^2-4)(r^2-9)\mu\delta^7
    \nonumber\\&&{}
    +\rat1{8!}r^2(r^2-1)(r^2-4)(r^2-9)\delta^8
    +\Ord{\delta^9}.
    \label{eq:doper}
\end{eqnarray}
The pattern in the above interpolation formula is clear.  
Now we
explore the numerical performance of a gap-tooth scheme using this
formula to determine teeth boundary conditions.

\begin{figure}
    \centering
    \includegraphics{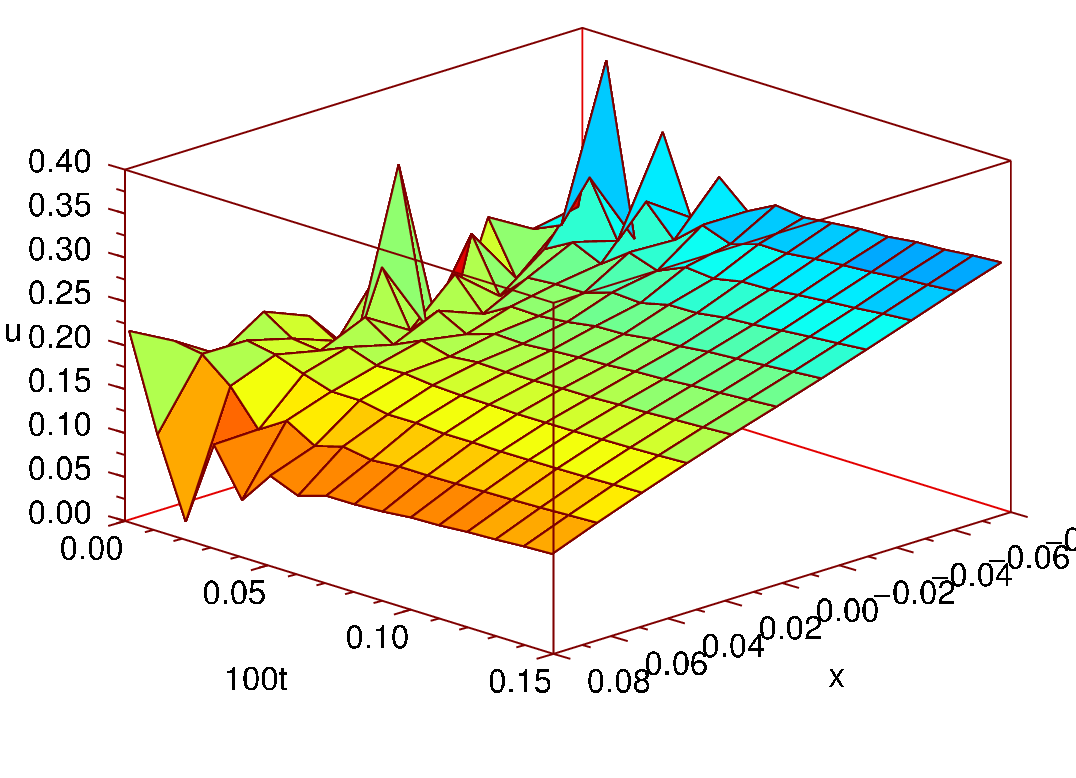}
	\caption{view of the initial microscopic evolution within a tooth
	with dynamics described by the diffusion \pde~\eqref{eq:diffn} and
	coupled to its neighbours.}
    \label{fig:b3z3}
\end{figure}

\begin{figure}
    \centering
    \includegraphics{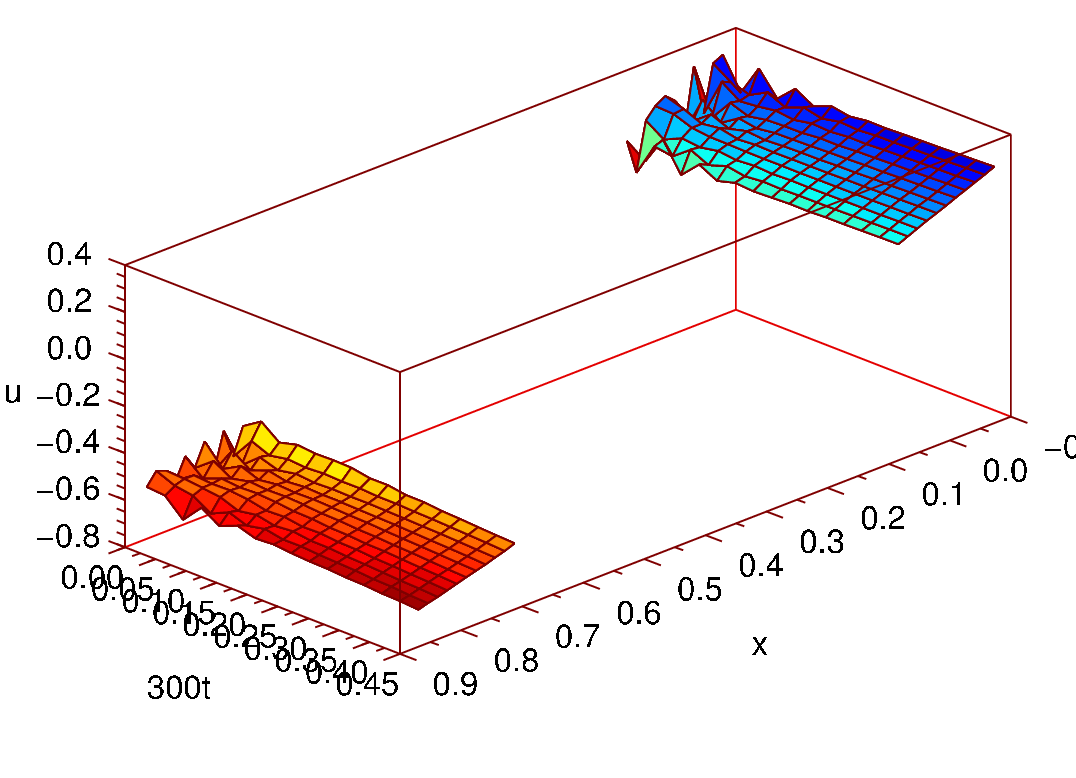}
	\caption{view of the initial microscopic evolution within a pair of
	neighbouring teeth with dynamics described by the diffusion
	\pde~\eqref{eq:diffn} and also coupled to their neighbours.}
    \label{fig:b3z4}
\end{figure}

Consider gap-tooth simulations of the simple
diffusion equation
\begin{equation}
    \D tu=\DD xu\,,
    \quad\text{and $2\pi$-periodic in~$x$.}
    \label{eq:diffn}
\end{equation}
Imagine we only have access to the dynamics through a microscopic
simulator of the diffusion~(\ref{eq:diffn}), here coded by a fine
discretisation on $n$~grid points, spaced a distance $\eta=h/(n-1)$
apart, across a tooth of microscopic width~$h=rH$\,.
The time integration is an explicit scheme with a 
microscopic time step, typically $\Delta t=10^{-6}$--$10^{-4}$.
Figure~\ref{fig:b3z3} shows an example of the initially rapid microscopic
evolution within one tooth; the microsimulator, coupled to its neighbors,
rapidly evolves to a smooth state.  
Figure~\ref{fig:b3z4}
similarly shows the initial evolution in two neighbouring teeth and
how the smooth subgrid field arises through the coupling to the
neighbouring teeth.  
Similar dynamics takes place during the initial
instants of the Burgers' evolution shown in Figure~\ref{fig:burg3}.

\begin{table}
    \centering
	\caption{Growth rates~$\lambda$ of perturbations from steady state
	$u=0$\,: for diffusion~(\ref{eq:diffn}) with $m$~teeth, $H=2\pi/m$\,;
	with gap to tooth ratio $r=0.1$\,; $n=11$ points in the microscale
	grid; and with the fourth order \pbc~(\ref{eq:pbc4}).}
    \label{tbl:spatb}
    \begin{tabular}{|r|llll|c|}
        \hline
        $m$ & \quad 1 & \quad 2,3 & \quad 4,5 & \quad 6,7 & $m+1:2m$  \\
        \hline
        4 & \E6{-12} & $-0.946256$ & $-2.166285$ & n/a & $-397.2$  \\
        8 & \E{-3}{-12} & $-0.996073$&$ -3.785024$&$ -7.121435$&$
        -1588.$ \\
        16 & \E{-1}{-10} & $-0.999750$&$ -3.984293$&$ -8.832102$&$
        -6355.$  \\
        32 & 0 & $-0.999986$&$ -3.998999$&$ -8.988613$&$ -25421.$  \\
        \hline
    \end{tabular}
\end{table}

Firstly we implement the following \pbc.
On the edge of the $j$th~tooth, at $x=X_j\pm rH$\,, the boundary condition
of the fine discretisation is that the 
field 
\begin{equation}
    v_j=\left[1 \pm r\mu\delta +\rat12r^2\delta^2 
    \pm \rat1{6}r(r^2-1)\mu\delta^3
    +\rat1{24} r^2(r^2-1)\delta^4 \right]U_j
    \,. \label{eq:pbc4}
\end{equation}
The first few terms of~(\ref{eq:doper}) provide this by interpolation from
the surrounding coarse grid values.
For the $j$th~tooth this \pbc\ involves macroscopic grid
values~$U_{j-2},\ldots,U_{j+2}$ only, and thus we should be able to
achieve $\Ord{H^4}$~consistency with the microsimulator.
We numerically linearize the map over one microscopic time step
by systematically perturbing each and every microscopic value from
zero (there are $mn$~such microscopic values, one for each of
$n$~fine grid points in each of $m$~teeth).
We then transform the eigenvalues~$\mu$ of this map to growth rates
$\lambda=\log(\mu)/\Delta t$\,.
The $mn$~growth rates fall into $n$~groups of $m$~modes.
Each group corresponds to a microscopic \emph{internal mode} of the dynamics;
the mode is essentially the same in each tooth.
Large negative growth rates correspond to rapidly decaying internal
modes with significant microscopic structure within each tooth.
The group of $m$~modes with \emph{small} growth rates correspond to the
relatively slowly evolving \emph{macroscopic} modes of interest that arise
through the coupling of the microscopic dynamics across the teeth.
Table~\ref{tbl:spatb} shows the leading seven growth rates, and the
magnitude of the leading internal growth rate, for various numbers of
teeth, $m=4,8,16,32$\,.
The exact growth rates of the diffusion \pde~(\ref{eq:diffn}) are
$\lambda=-k^2$ for integer~$k$.
The table shows that as the number of teeth doubles, the accuracy of the
growth rates of the macroscopic modes improves by a factor of
about~$16$.
This is consistent with an~$\Ord{H^4}$ method as predicted for
diffusion with \pbc~(\ref{eq:pbc4}).

\begin{table}
    \centering
	\caption{Growth rates~$\lambda$ of perturbations from steady state
	$u=0$\,: for diffusion~(\ref{eq:diffn}) with $m$~teeth, $H=2\pi/m$\,;
	with gap to tooth ratio $r=0.1$\,; $n=11$ points in the microscale
	grid; and with the sixth order \pbc.}
    \label{tbl:spatc}
    \begin{tabular}{|r|llll|c|}
        \hline
        $m$ & \quad 1 & \quad 2,3 & \quad 4,5 & \quad 6,7 & $m+1:2m$  \\
        \hline
        4 & \E{-5}{-12} & $-0.981981$ & $-2.453767$ & n/a & $-397.2$  \\
        8 & \E{1}{-11} & $ -0.999653$&$ -3.927925$&$ -7.835158$&$ -1588.$ \\
        16 & \E{8}{-11} & $-1.000001$&$ -3.998611$&$ -8.966332$&$ -6355.$  \\
        32 & \E{8}{-10} & $-1.000002$&$ -4.000004$&$ -8.999518$&$ -25421.$  \\
        \hline
    \end{tabular}
\end{table}

Table~\ref{tbl:spatc} shows the even higher order accuracy from
implementing sixth order \pbc{}s from~(\ref{eq:doper})---growth
rates slightly larger than the ideal seem to be due to the relatively
small number of microscopic grid points within the teeth.  
These sixth order \pbc{}s are
used in the simulations of the nonlinear Burgers' equation shown in
Figure~\ref{fig:burg3}.  
This simulation suggests that gap-tooth schemes
employing such \pbc{}s even for nonlinear systems are 
effective.

\section{Mixed boundary conditions for the teeth}
\label{sec:mixed}

Let us explore mixed boundary conditions at the edges of the teeth:
suppose the microsimulator requires $av_j\pm b\partial_xv_j$ specified
on the edge of the teeth $x=x_j\pm rH$ for some constants $a$~and~$b$.
The case $a=1$ and $b=0$ constitutes Dirichlet \pbc{}s discussed in the
previous section.
The case $a=0$ and $b=1$ constitutes Neumann \pbc{}s as discussed by
\cite{Roberts04d}: there we used the interpolation
formula~\eqref{eq:boper} to specify slopes\slash fluxes on the edge of each tooth;
we obtained spectra of accuracy similar to those in Tables
\ref{tbl:spatb}~and~\ref{tbl:spatc}.

For mixed \pbc{}s we propose to simply combine
(\ref{eq:boper})~and~(\ref{eq:doper}) to give, for example, the fourth
order in macroscopic grid size~$H=2\pi/m$ boundary condition   
\begin{eqnarray}&&
    av_j\pm b\partial_xv_j
    \nonumber\\&&{}
    =a\left[ 1 \pm r\mu\delta +\rat12r^2\delta^2 
    \pm \rat1{6}r(r^2-1)\mu\delta^3
    +\rat1{24} r^2(r^2-1)\delta^4 \right] U_j
    \nonumber\\&&{}
    \pm \frac{b}H\left[ \mu\delta \pm r\delta^2 
    -(\rat16-\rat12r^2)\mu\delta^3
    \mp r(\rat1{12}-\rat16r^2)\delta^4 \right]U_j
\nonumber \\ && \text{on } x=x_j\pm rH\,.
    \label{eq:mbc4}
\end{eqnarray}
\begin{table}
    \centering
	\caption{Growth rates~$\lambda$ of perturbations from steady state
	$u=0$\,: for diffusion~(\ref{eq:diffn}) with $m$~teeth, $H=2\pi/m$\,;
	with gap to tooth ratio $r=0.1$\,; $n=11$ points in the microscale
	grid; and with the mixed \pbc~(\ref{eq:mbc4}) with
	$a=0.95$ and $b=0.05$\,.}
    \label{tbl:spatm}
    \begin{tabular}{|r|llll|c|}
        \hline
        $m$ & \quad 1 & \quad 2,3 & \quad 4,5 & \quad 6,7 & $m+1:2m$  \\
        \hline
        4 & \E9{-12} & $-0.939448$ & $-2.151772$ & n/a &
        $-240.7$ \\
        8 & \E{-2}{-11} & $-0.990854$ & $-3.766007$ & $-7.089004$ &
        $-756.9$ \\
        16 & \E7{-11} & $-0.996405$ & $-3.971027$ & $-8.803476$ & 
        $-2417.$ \\
        32 & \E6{-10} & $-0.998047$ & $-3.991243$ & $-8.971213$ &
        $-8153.$ \\
        \hline
    \end{tabular}
\end{table}%
We use $a=0.95$ and $b=0.05$ in the mixed \pbc: this gives a mixed
boundary condition where the effects of the function value~$v_j$ and
its gradient~$\partial_xv_j$ are roughly comparable in the \pbc\ (if the
parameter~$b$ is significantly larger, then the gradient term dominates
the \pbc).
The numerical eigenvalues given in Table~\ref{tbl:spatm} for the
diffusion equation~(\ref{eq:diffn}) with these \pbc{}s again show
convergence to the correct eigenvalues as the number of teeth increases,
that is, as the macroscopic grid size~$H\to0$\,.
However, the convergence is not as rapid as for Dirichlet \pbc{}s.
The poorer convergence as $H\to0$ seems to be due to the microscopic
grid resolution: successive doubling of the number of interior points,
see Table~\ref{tbl:spatmm}, demonstrates that there are significant errors
of~$\Ord{\eta^2}$ in the microscale grid size~$\eta$.  
Thus the total
error in this implementation of the mixed \pbc{}s seems to
be~$\Ord{H^4,\eta^2}$.

\begin{table}
    \centering
	\caption{Growth rates~$\lambda$ of perturbations from steady state
	$u=0$\,: for diffusion~(\ref{eq:diffn}) with $m=8$~teeth; with gap
	to tooth ratio $r=0.1$\,; $n$~points in the microscale grid to
    show variation with microscale resolution;
    and with the mixed \pbc~(\ref{eq:mbc4}) with
    $a=0.95$ and $b=0.05$\,.}
    \label{tbl:spatmm}
    \begin{tabular}{|r|llll|}
        \hline
        $n$ & \quad 1 & \quad 2,3 & \quad 4,5 & \quad 6,7   \\
        \hline
        11 & \E{-2}{-11} & $-0.990854$ & $-3.766007$ & $-7.089004$ \\
        21 & \E1{-10} & $-0.994896$ & $-3.781289$ & $-7.117588$ \\
        41 & \E{-3}{-10} & $-0.995792$ & $-3.784677$ & $-7.123924$ \\
        \hline
    \end{tabular}
\end{table}

Here the microscale simulation is that of a fine discretisation of a \pde.
Thus the derivatives in the mixed \pbc~\eqref{eq:mbc4} are subject to
the significant errors of numerical differentiation when computed on
the microscale.
%
%
As Table~\ref{tbl:spatmm} shows, the approximation of
derivatives does incur errors; we would be better off without such
errors.
Higher order formulae for microscale interpolation would reduce the
microscale errors in the boundary derivatives, perhaps from
$\Ord{\eta^2}$~to~$\Ord{\eta^4}$, but would ruin the small
bandwidth of the microscale simulator.
In any case, recall that we adopt the policy that we cannot change the
microscale simulator as it is a legacy code handed to us from past
development.
We cannot (do not want to) change the nature nor accuracy of its boundary
conditions.
Consequently we proceed to address the problem of supplying boundary
conditions at the edge of the teeth, \emph{precisely as required 
during execution} of the legacy microscale simulator.

\section{Teeth with two point boundary conditions}
\label{sec:twopt}

A microscale simulator may have implemented boundary conditions that do
not fit into the classic partial differential equation form of
Dirichlet, Neumann nor mixed.
Here the microscale simulator implements a discretisation of 3~point stencil
width.
Consequently the simulator has been written so that the supplied
boundary conditions only depend upon each of the two extreme pairs of
points in each tooth.
We thus investigate teeth boundary conditions that specify a combination
of these two point values of the field at the edge of each tooth.
This specific case is just one example of the wide range of possible
nonlocal \pbc{}s that specific microsimulators may require.

Suppose the microsimulator, here a fine spatial discretisation of the
diffusion \pde~\eqref{eq:diffn}, implements a tooth boundary
condition of the (linear) form
\begin{equation}
     v_{j,1}+\beta v_{j,2}
    \quad\text{and}\quad
    \beta v_{j,n-1}+ v_{j,n}
    \quad\text{are specified,}
    \label{eq:bcg}
\end{equation}
where $v_{j,i}$ denotes the microscale field value at the $i$th~microscale
grid point in the $j$th~tooth.
For example, the case $\beta=1$ approximates Dirichlet boundary
conditions at the microscale grid mid-points $x_{j,3/2}$
and~$x_{j,n-1/2}$ near the edges of each tooth, which is exactly the
case $a=0$ and $b=1$ implemented in the previous
Section~\ref{sec:mixed}.
Different values of~$\beta$ would \emph{approximate} different mixed boundary
conditions of the previous section.

The procedure is straightforward: we interpolate the macroscale grid
values to find the specific values required by the boundary
conditions~(\ref{eq:bcg}).
Recall that~(\ref{eq:pbc4}) gives a fourth order interpolation from the
macroscale grid to points at the tooth boundaries $x=X_j\pm rH$\,; this
gives appropriate values for $v_{j,1}$~and~$v_{j,n}$.
Get appropriate values for $v_{j,2}$~and~$v_{j,n-1}$ through simply
replacing in the formula the ratio~$r=h/(2H)$ by the ratio required to reach the penultimate
microgrid point, namely $r'=(h/2-\eta)/H$, where $\eta=h/(n-1)$ is
the microgrid size.
Thus the fourth order version of the boundary condition~(\ref{eq:bcg})
is that at $x=X_j\pm rH$
\begin{eqnarray}&&
    (1+ \beta E^{\mp\eta/H})v_j
    \nonumber\\
    &&= \left[1 \pm r\mu\delta +\rat12r^2\delta^2 
    \pm \rat1{6}r(r^2-1)\mu\delta^3
    +\rat1{24} r^2(r^2-1)\delta^4 \right]U_j
    \label{eq:pbg4}\\&&{}
    +\beta \left[1 \pm {r'}\mu\delta +\rat12{r'}^2\delta^2 
    \pm \rat1{6}{r'}({r'}^2-1)\mu\delta^3
    +\rat1{24} {r'}^2({r'}^2-1)\delta^4 \right]U_j
    \,. \nonumber
\end{eqnarray}
\begin{table}
    \centering
	\caption{Growth rates~$\lambda$ of perturbations from steady state
	$u=0$\,: for diffusion~(\ref{eq:diffn}) with $m=8$~teeth; with gap
	to tooth ratio $r=0.1$\,; $n$~points in the microscale grid to
    show variation with microscale resolution $\eta\propto1/n$\,;
    and with the fourth order general \pbc~(\ref{eq:pbg4}) with
    $\beta=1$\,.}
    \label{tbl:spatgg}
    \begin{tabular}{|r|llll|}
        \hline
        $n$ & \quad 1 & \quad 2,3 & \quad 4,5 & \quad 6,7   \\
        \hline
        11 & \E1{-10} & $-0.999741$ & $-3.984137$ & $-8.831209$ \\
        21 & \E{-3}{-10} & $-0.999742$ & $-3.984159$ & $-8.831368$ \\
        41 & \E1{-9} & $-0.999742$ & $-3.984169$ & $-8.831453$ \\
        \hline
    \end{tabular}
\end{table}%
Implementing the \pbc~(\ref{eq:pbg4}) for the diffusion
equation~(\ref{eq:diffn}) gives a numerical approximation scheme with
eigenvalues shown in Table~\ref{tbl:spatgg} for varying microgrid
resolution.
See that there is only an extremely weak dependence upon the microgrid
size~$\eta$.
\emph{Thus implementing directly the boundary conditions that the
microscale simulator actually expects during execution results in 
much better accuracy than
trying to approximate the microscale \pbc{}s using computed spatial
derivatives.}

\begin{table}
    \centering
	\caption{Growth rates~$\lambda$ of perturbations from steady state
	$u=0$\,: for diffusion~(\ref{eq:diffn}) with $m$~teeth,
	$H=2\pi/m$\,; with gap to tooth ratio $r=0.1$\,; $n=11$ points in
	the microscale grid; and with the fourth order general
	\pbc~(\ref{eq:pbg4}) with $\beta=1$\,.}
    \label{tbl:spatg}
    \begin{tabular}{|r|llll|c|}
        \hline
        $m$ & \quad 1 & \quad 2,3 & \quad 4,5 & \quad 6,7 & $m+1:2m$  \\
        \hline
        4 & \E{-8}{-12} & $-0.946069$ & $-2.165068$ & n/a & $-489.5$ \\
        8 & \E4{-11} & $-0.996034$ & $-3.784277$ & $-7.118312$ & $-1958.$ \\
        16 & \E1{-10} & $-0.999741$ & $-3.984137$ & $-8.831209$ & $-7832.$ \\
        32 & \E8{-10} & $-0.999983$ & $-3.998964$ & $-8.988427$ & $-31329.$ \\
        \hline
    \end{tabular}
\end{table}%

Lastly, Table~\ref{tbl:spatg} shows the eigenvalues of the gap-tooth
scheme for varying number~$m$ of teeth in the domain.
See that the eigenvalues converge to their correct values
like~$\Ord{H^4}$ as expected by the construction.

Higher order \pbc{}s, in the macroscopic grid size~$H$, would similarly be
based upon the expansion~(\ref{eq:doper}).
We then expect even more rapid convergence as the macroscale grid size $H\to
0$\,.

\section{The model is independent of the tooth boundary conditions}
\label{sec:ind}

Here we use analytic methods of holistic discretisation
\cite[e.g.]{Roberts01g} to explore the gap-tooth scheme on a general
class of \pde{}s with general mixed boundary conditions.
The analysis establishes three important properties:
\begin{itemize}
    \item the approach generates macroscopic models which are
    consistent with the microscopic dynamics to high orders in grid
    spacing~$H$;
    
	\item the macroscopic model and the microscopic solution field are
	essentially independent of the size of the teeth, as parametrised
	by~$r$; and

	\item the macroscopic model and the microscopic solution field are
	essentially independent of the details of the \pbc{}s.
\end{itemize}

\subsection{Theory underpins analysis of a PDE with tooth boundary
conditions}

We explore solutions of the class of \emph{linear}
hyper-advection-diffusion \pde{}s
\begin{equation}
    \D tu=\DD xu-\epsilon \left(c\D xu+b\Dn x3u+a\Dn x4u \right)\,,
    \label{eq:pdeg}
\end{equation}
where $a$, $b$~and~$c$ are arbitrary parameters, and where
$\epsilon$~is introduced as a convenient mechanism to control
truncation in the multivariate power series solutions in the
parameters~$a$, $b$~and~$c$.
This \pde\ is solved with mixed tooth boundary conditions inspired
by~(\ref{eq:mbc4}), namely that on $x=x_j\pm rH$\,, and in terms of an
artificial parameter~$\gamma$ that we explain shortly,
\begin{eqnarray}&&
    \pm \alpha v_j+\partial_xv_j
    \nonumber\\&&{}
    =\pm\alpha\left\{ 1 
    +\gamma r\left[ \pm \mu\delta +\rat12r\delta^2 \right]
    +\gamma^2\rat1{6}r(r^2-1)\left[ \pm \mu\delta^3
    +\rat1{4} r\delta^4 \right]
    \right.\nonumber\\&&\quad\left.{}
    +\gamma^3\rat1{120}r(r^2-1)(r^2-4)\left[ \pm \mu\delta^5
    +\rat1{6} r\delta^6 \right]
    \right\} U_j
    \nonumber\\&&{}
    + \frac{1}H\left\{ \gamma\left[ \mu\delta \pm r\delta^2 \right]
    +\gamma^2\left[ -(\rat16-\rat12r^2)\mu\delta^3
    \mp r(\rat1{12}-\rat16r^2)\delta^4 \right]
    \right.\nonumber\\&&\quad\left.{}
    +\gamma^3\left[ (\rat1{30}-\rat18r^2+\rat1{24}r^4)\mu\delta^5
    \pm r(\rat1{90}-\rat1{36}r^2+\rat1{120}r^4)\delta^6 \right] \right\}U_j
    \nonumber\\&&
    +\Ord{\gamma^4}\,.
    \label{eq:mbc}
\end{eqnarray}
Explore the structure of this complicated looking \pbc: $\pm \alpha
v_j+\partial_xv_j$ represents a general linear combination of the
microscopic field at the edge of each tooth that needs to be specified
for the microscopic simulator; those terms in the right-hand side
multiplied by~$\pm\alpha$ form the estimate of the field~$v_j$
interpolated from the surrounding macroscopic grid values; those terms
in the right-hand side multiplied by~$1/H$ form the estimate of the
field's gradient~$\partial_xv_j$ interpolated from the surrounding
macroscopic grid values.
However, these two interpolations only hold when the artificial
parameter $\gamma=1$\,; one is the physically interesting
value of~$\gamma$.
Why then do we introduce the parameter~$\gamma$? 
The reason is that, as
in ``discretisation'' \cite[]{Roberts98a}, based around the
special values of the parameters $\gamma=\alpha=\epsilon=0$\,, the
general \pde~(\ref{eq:pdeg}) with \pbc~(\ref{eq:mbc}) possesses a (slow)
centre manifold parametrised by the macroscopic grid values.
On this centre manifold the evolution of these macroscopic grid values
forms a macroscale model of the \pde.
This model has rigorous theoretical support based upon
$\gamma=0$\,, and it becomes physically relevant when evaluated at $\gamma=1$.

We briefly explain how centre manifold theory underpins the macroscale
model.
Initially set $\gamma=\alpha=\epsilon=0$; then the \pde+\pbc\ become the
diffusion equation with \emph{insulating} boundaries at the edges of the
teeth, $x=x_j\pm rH$\,.
Thus, exponentially quickly, all structure within each tooth diffuses
away to become constant, but a different constant for each tooth depending
upon the initial conditions.
See a similar evolution in Figures \ref{fig:b3z3}~and~\ref{fig:b3z4};
but there the teeth are coupled, so that the rapid evolution is to a smooth
variation in each tooth, whereas here the insulated evolution,
$\gamma=\alpha=\epsilon=0$\,, is to a constant in each tooth.
%
%
But we are only interested in fully coupled teeth for which
$\gamma=1$\,, and in non-zero $\alpha$~and~$\epsilon$.
Thus from the simple base of piecewise constant fields, we construct a
description of the field~$u$ and its slow evolution as a power series
in the ``perturbations'' measured by $\gamma$, $\alpha$~and~$\epsilon$.
The departure of the field~$u$ from a constant within each tooth gives
the microscopic (subgrid, subtooth) field, as shown for example in the smooth fields of
Figures \ref{fig:b3z3}~and~\ref{fig:b3z4} that are quickly established.
The slow evolution of the coarse grid values~$U_j$ gives the macroscopic
model.

The various powers of~$\gamma$ in the \pbc~(\ref{eq:mbc}) are chosen so
that truncation of the expressions to errors~$\Ord{\gamma^p}$ will
generate a discrete macroscopic model expressing $\dot U_j$ in terms of
only $U_{j-p+1},\ldots,U_{j+p-1}$ (a spatial stencil of width~$2p-1$).
Centre manifold theory \cite[e.g.]{Carr81, Kuznetsov95} asserts that
\begin{itemize}
    \item  such a model exists,

	\item that through its exponential attractiveness, the model is
	relevant in some finite neighbourhood of
	$\gamma=\alpha=\epsilon=0$\,,

	\item and that we may systematically construct the power series
	approximation to the model.
\end{itemize}
Because truncation to errors~$\Ord{\gamma^p}$ results in a model with
stencil width~$2p-1$, \emph{such a truncation corresponds to the
gap-tooth scheme utilising \pbc{}s involving interpolation from only
the  $2p-1$ neighbouring grid values $U_{j-p+1},\ldots,U_{j+p-1}$.}

Computer algebra\footnote{\protect \url
{http://www.sci.usq.edu.au/staff/aroberts/CA/burgermixed.red} is
the source script which was available at the time of writing.} performs
all the tedious details of constructing the model \cite[]{Roberts96a}.
We seek a model where the subtooth\slash subgrid field 
\begin{equation}
    u(x,t)=v_j(x,\vec U;\alpha,\gamma,\epsilon)\,.
    \label{eq:cmm}
\end{equation}
That is, the subtooth field has some spatial structure, such as that in
Figures \ref{fig:b3z3}~and~\ref{fig:b3z4}, which: depends upon the
neighbouring grid values $U_{j-p+1},\ldots,U_{j+p-1}$\,; may depend
upon the specific \pbc{} through its parameter~$\alpha$; depends upon
the specific \pde{} though its parameter~$\epsilon$; and depends upon
the coupling parameter~$\gamma$.  
Centre manifold theory assures us the
evolution of the system is governed by the evolution of the grid
values:
\begin{equation}
    \dot U_j=g_j(\vec U;\alpha,\gamma,\epsilon)\,;
    \label{eq:cmu}
\end{equation}
this formula is the macroscopic (closed) discretisation.  
We solve by iteration
the \pde~(\ref{eq:pdeg}) with \pbc~(\ref{eq:mbc}) to find the centre
manifold~\eqref{eq:cmm} and its associated
coarse discretisation~\eqref{eq:cmu}.  
The results are expressions for the
microscopic fields~$v_j$ and the macroscopic evolution, $\dot
U_j=g_j$\,, that are accurate to some specified order in the small
parameters $\alpha$, $\gamma$~and~$\epsilon$.

\subsection{Modelling $\Ord{\epsilon}$ changes to the PDE}

The subgrid fields in each tooth will in general depend upon the
coefficients $a$, $b$~and~$c$ that determine the \pde.
For example, there are nontrivial dependencies upon the advection speed
that ensure the macroscale model naturally transforms to an upwind
discretisation for large advection speeds~$c$
\cite[]{Roberts00a}---such influences show up in the
$\Ord{\epsilon^2}$~terms that we explore in the next subsection.
Here we first explore the models linear in $a$,
$b$~and~$c$, that is, linear in the general changes to the
\pde~\eqref{eq:pdeg} with the \pbc~(\ref{eq:mbc}).

For example, to errors $\Ord{\alpha^3,\gamma^4,\epsilon^2}$,
computer algebra generates the macroscopic evolution 
\begin{eqnarray}
    \dot U_j &=&
    \frac1{H^2}\left[ \gamma\delta^2 -\rat1{12}\gamma^2\delta^4
    +\rat1{90}\gamma^3\delta^6 \right] U_j
    \nonumber\\&&{}
    -\frac{\epsilon c}{H}\left[ \gamma\mu\delta 
    -\rat16\gamma^2\mu\delta^3 +\rat1{30}\gamma^3\mu\delta^5 \right] U_j
    \nonumber\\&&{}
    -\frac{\epsilon b}{H^3}\left[ \gamma^2\mu\delta^3 
    -\rat14\gamma^3\mu\delta^5 \right] U_j
    \nonumber\\&&{}
    -\frac{\epsilon a}{H^4}\left[ \gamma^2\delta^4
    -\rat16\gamma^3\delta^6 \right] U_j
    +\Ord{\alpha^3,\gamma^4,\epsilon^2}\,.
    \label{eq:lowm}
\end{eqnarray}
When evaluated at the physically relevant parameter $\gamma=1$ these
are the classical finite difference operators for the
\pde~(\ref{eq:pdeg}), truncated to~$\Ord{\delta^7}$.
Consequently the terms in the macroscopic model~\eqref{eq:lowm} are
\emph{consistent} with the \pde~(\ref{eq:pdeg}) to various orders in
the macroscopic grid size~$H$.
The order of consistency depends upon the order of truncation in the
artificial coupling parameter~$\gamma$ and the order of the derivatives
in each term.
\emph{Observe that the macroscopic evolution operator is
\emph{independent} of~$r$, the size of the teeth, and
\emph{independent} of~$\alpha$ which parametrises the precise nature of
the \pbc~(\ref{eq:mbc}).}

Now we explore the microscopic field within the teeth.
To low order in the coupling parameter~$\gamma$ and in terms of the
microscopic tooth space variable $\xi=(x-X_j)/H$\,, we find
\begin{eqnarray}
    v_j&=&U_j +\gamma\left\{
    \left[ \xi\mu\delta +\half\xi^2\delta^2 \right]
    \right.\nonumber\\&&\left.{}
    +\epsilon c H\left[ (\rat16\xi^3-\half r^2\xi)
    +\rat13 H\alpha r^3\xi 
    -\rat13 H^2\alpha^2r^4\xi \right]\delta^2
    \right\}U_j 
    \nonumber\\&&{}
    +\Ord{\alpha^3,\gamma^2,\epsilon^2}\,.
    \label{eq:mlo}
\end{eqnarray}
The first line gives the classic quadratic interpolation through the
grid values $U_j$~and~$U_{j\pm1}$.
The second line shows microscopic field structure in the advection
speed~$c$.
But it exhibits undesirable dependence upon the tooth width~$r$ and
nature~$\alpha$ of the \pbc.
However, inspect the next order terms in coupling parameter~$\gamma$:
\begin{eqnarray}
    v_j&=&\cdots +\gamma^2\left\{
    \left[ \rat16(-\xi+\xi^3)\mu\delta^3 
    +\rat1{24}(-\xi^2+\xi^4)\delta^4 \right]
    \right.\nonumber\\&&\left.{}
    +\epsilon c H\left[ -(\rat16\xi^3-\half r^2\xi)
    -\rat13 H\alpha r^3\xi
    +\rat13 H^2\alpha^2r^4\xi \right]\delta^2
    \right.\nonumber\\&&\left.{}
    +(-\rat1{18}\xi^3 +\rat1{60}\xi^5 +\rat16r^2\xi
    -\rat1{12}r^2\xi^3 \rat16r^4\xi )\delta^4
    \right.\nonumber\\&&\left.{}
    +\alpha H(-\rat19r^3\xi +\rat1{18}r^3\xi^3 -\rat4{15}r^5\xi )\delta^4
    \right.\nonumber\\&&\left.{}
    +\alpha^2 H^2(\rat19r^4\xi -\rat1{18}r^4\xi^3 +\rat{17}{45}r^6\xi
    )\delta^4
    \right.\nonumber\\&&\left.{}
    +\frac{\epsilon b} H\left[ -(\rat16\xi^3-\half r^2\xi)
    -\rat13 H\alpha r^3\xi
    +\rat13 H^2\alpha^2r^4\xi \right]\delta^4
    \right\}U_j 
    \nonumber\\&&{}
    +\Ord{\alpha^3,\gamma^3,\epsilon^2}\,;
    \label{eq:medm}
\end{eqnarray}
The dots denote the terms given in the right-hand side
of~(\ref{eq:mlo}).
The first line in the above higher order terms contains reassuringly the
classic quartic interpolation formulae.
The second line, when we set $\gamma=1$\,, cancels all the undesirable
$\alpha$~and~$r$ dependence in the lower order~(\ref{eq:mlo}).
The third and later lines above describe higher order microscopic
structure; some of this undesirably depends upon the \pbc{}
through~$\alpha$ and the size of the teeth through~$r$, but as far as
we have explored, any $\alpha$~and~$r$ dependence introduced at any
order in~$\gamma$ is canceled by terms at higher orders in~$\gamma$.
Thus any finite truncation of the power series expansion for the model
may have undesirable dependence upon $\alpha$~and~$r$, but as the order
in the artificial parameter~$\gamma$ is increased, this dependence is
removed.
In this sense, the microscopic field is ``essentially'' independent of
the details of the \pbc\ and independent of the tooth width~$r$.

\subsection{Modelling~$\Ord{\epsilon^2}$ effects in the PDE}

In our analysis we find that even order operators in the microscale
\pde, such as the diffusion~$u_{xx}$ and the
hyper-diffusion~$au_{xxxx}$, are represented simply in the macroscale
discretisation.
However, odd order operators in the microscale \pde, such as
advection~$cu_x$ and the dispersion~$bu_{xxx}$, create nontrivial
effects; these first show up in terms quadratic in their amplitude and
hence they first appear in terms of~$\Ord{\epsilon^2}$.
We now show that the effects of $\Ord{\epsilon^2}$~terms typically act
to stabilise the $\Ord{\epsilon}$ discrete model.
The implication for the gap-tooth method is that the resolution of the
subgrid structures by the microscale simulator will also typically
maintain stability of the discrete macroscale model.
That is, the microscale simulator will provide successful closure for
the macroscale discretisation when coupled with the proposed \pbc{}s.


To illustrate this, we can construct the approximate model of the
\pde~(\ref{eq:pdeg}) with \pbc~(\ref{eq:mbc}) to errors
$\Ord{\alpha^2,\gamma^4,\epsilon^3}$; that is, we include quadratic
effects in the coefficients $a$, $b$~and~$c$.
The details of the model are too long to record here.
However, we find that the equivalent \pde\ to the macroscale discrete
model~(\ref{eq:lowm}) with its $\Ord{\epsilon^2}$ modifications is
\begin{eqnarray}
    &&\D tu =
    \gamma \DD xu -\epsilon\left[\gamma c\D xu +\gamma^2 b\Dn x3u 
    +\gamma^2a\Dn x4u \right]
    \nonumber\\&&{}
    +\epsilon^2\left[ -(\gamma-\gamma^2)bc\DD xu
    -(\gamma^2-\gamma^3)b^2\Dn x4u \right]
    \nonumber\\&&{}
    +H^2\left\{ \rat1{12}(\gamma-\gamma^2)\Dn x4u
    \right.\nonumber\\&&\left.\quad{}
    +\epsilon\left[ -\rat16(\gamma-\gamma^2)c\Dn x3u
    -\rat14(\gamma^2-\gamma^3)b\Dn x5u
    -\rat16(\gamma^2-\gamma^3)a\Dn x6u \right]
    \right.\nonumber\\&&\left.\quad{}
    +\epsilon^2\left[ +\rat13(\gamma-\gamma^2)c^2r^2\DD xu
    -\rat1{12}(\gamma-5\gamma^2+4\gamma^3)bc\Dn x4u
    \right.\right.\nonumber\\&&\left.\left.\quad\quad{}
    +\rat23(\gamma^2-\gamma^3)bcr^2\Dn x4u
    -\rat16(\gamma^2 -7\gamma^3+4\gamma^3r^2)b^2\Dn x6u \right] \right\}
    \nonumber\\&&{}
    +\Ord{H^3,\alpha^2,\gamma^4,\epsilon^3}\,.
    \label{eq:quad}
\end{eqnarray}
Observe that to this level of accuracy there is no dependence upon the
\pbc\ parameter~$\alpha$, thus our comments apply for all the \pbc{}s.
Now consider the components of~(\ref{eq:quad}) in turn.
The first line of~(\ref{eq:quad}) is the original general linear
\pde~(\ref{eq:pdeg}) when evaluated at the physically meaningful
$\gamma=1$\,.
The second line shows some error terms, quadratic in~$\epsilon$, that
disappear for $\gamma=1$\,.
The $bcu_{xx}$ error disappears when $\Ord{\gamma^2}$~terms are
retained, which is as soon as the discretisation stencil is wide
enough to model the third order dispersion term~$bu_{xxx}$.
The $b^2u_{xxxx}$ dispersion-induced term shows that the method initially
incorporates its effects as enhanced dissipation, as the coefficient
of~$\gamma^2$ is negative; then, when higher order accuracy is
requested by retaining $\Ord{\gamma^3}$~terms, the method proceeds to remove the
incurred error via the $\gamma^3$~term.

Consider applying these \pbc\ to a microscale simulator; 
for slow enough spatial variations,
the microsimulator is equivalent to some `infinite order' \pde.
For example, the \emph{microscale} discretisation $\dot u_i=(1/\eta^2)
\delta^2 u_i$ is, by~\eqref{eq:derdel}, equivalent to the \pde\ $u_t=
(4/\eta^2)\sinh^2(\eta\partial_x/2)u$\,.
We expect the behaviour of errors in the gap-tooth scheme seen here,
induced by the $u_{xxx}$~and~$u_{xxxx}$ terms,
be representative of the behaviour of errors in the `high order'
equivalent terms of any given microscopic simulator.

The remaining terms in~(\ref{eq:quad}) are $\Ord{H^2}$ and hence vanish
as the macroscopic grid size $H\to 0$\,.
Nonetheless look at the $H^2$ terms in the third to sixth lines as they
apply to simulations with finite~$H$.
The third and fourth lines show that the initial discretisation errors of
the linear terms are eliminated, for the physical $\gamma=1$\,, via
the next higher order in coupling parameter~$\gamma$.
The particular focus of this subsection is the $\epsilon^2$~terms on
the next two lines.
The $c^2u_{xx}$ terms show that at low order truncations in~$\gamma$
the method treats advection in a manner that increases dissipation, as
the coefficient is positive, and thus helps to maintain the stability of
the macroscale discretisation to high advection speeds~$c$
\cite[explored in][]{Roberts98a,Roberts00a}.
The interaction between advection~$cu_x$ and dispersion~$bu_{xxx}$ can
maintain stability or be destabilising depending upon the sign
of~$bc$: whether we truncate at~$\gamma$ when $-bcu_{xx}$ dominates
(second line) or we truncate at $\gamma^2$ when $bcu_{xxxx}$ dominates
(fifth and sixth lines), the combination is stabilising whenever
$bc<0$\,, that is, when the phase velocity of wave-like effects does
not change direction as a function of wavenumber.
The last term on the sixth line will be dominated by the dissipative
$b^2u_{xxxx}$ term on the second line, and we presume will vanish at
higher orders in the coupling parameter~$\gamma$.
Thus, from the equivalent \pde~(\ref{eq:quad}) of the macroscale model,
we deduce that the subgrid scale interactions between processes in the
\pde, and hence for microscale simulators in general, are accounted for
in this approach to generate a macroscale model that is typically
stable.

Indeed this equivalent \pde~(\ref{eq:quad}) confirms support for the
gap-tooth scheme with \pbc\ by centre manifold theory.  
Theory asserts that the original system, here the
\pde~\eqref{eq:quad}, and the centre manifold model, here the
macroscale discretisation~\eqref{eq:cmu}, have the same stability.
Thus when the microscale system is stable, so will the macroscale
discretization.  
The caveat is that we can only construct the centre
manifold approximately; we control the errors to some order in the
parameters $\alpha$ $\gamma$~and~$\epsilon$, but there will be some
error, albeit of high order in the parameters.

\section{Conclusion}


\begin{figure}[btp]
    \centering
    \includegraphics{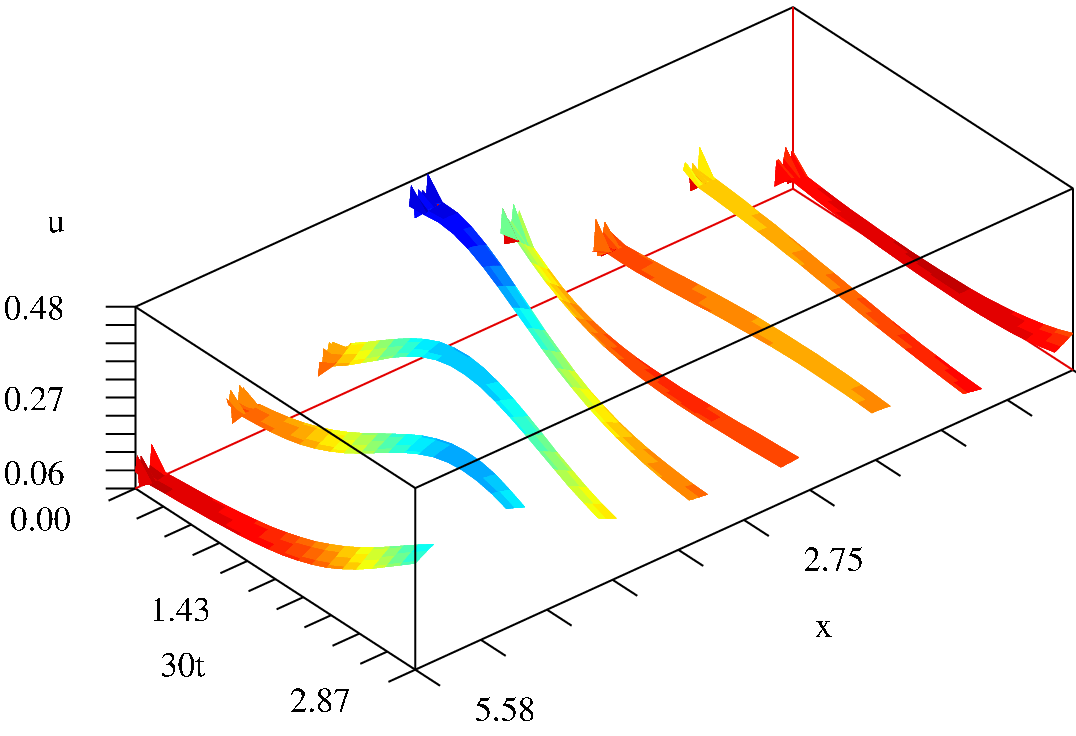}
	\caption{simulation of Burgers' equation using general 2-point
	boundary condition on the teeth of the fourth order~(\ref{eq:pbg4})
	with $\beta=1$ demonstrates the method is stable even for
    nonlinear \pde{}s.}
    \label{fig:burg3b}
\end{figure}

We use macroscale interpolation based upon the expansions
(\ref{eq:boper})~and~(\ref{eq:doper}) to determine \pbc{}s for the
boundary conditions at the edge of the teeth in the gap-tooth scheme.
The interpolation was used to implement \emph{directly} whatever boundary
conditions are \emph{actually} needed by the microscale legacy code
during execution.
Figure~\ref{fig:burg3b} shows a simulation of the nonlinear Burgers'
equation with 2~point boundary conditions at the boundaries of each
tooth as an illustrative example.
We found that the macroscopic models resulting from the microsimulator
and the constructed \pbc\ were consistent, to high order, with the microscopic
dynamics; that the macroscopic models and the microscopic (subgrid, subtooth)
fields were essentially independent of the tooth size and the detailed
nature of the \pbc.
We expect the same type of \pbc\ to be effective for microsimulations in
more than one spatial dimension.
Interesting future research would seek \pbc\
that do not require communication across the gaps 
between the teeth at each and every microscale time step,
and the interplay of \pbc with implicit integration schemes.

Further exciting research would explore issues of existence
and performance of \pbc{}s for stochastic microsimulators.


\paragraph{Acknowledgment:} I.~G.~K. is supported in part by
\textsc{darpa} and the \textsc{US DOE}.

\bibliographystyle{agsm}
\bibliography{news,ajr,bib,new}

\end{document}